\documentclass[12pt,a4paper,twoside]{article}

\newcommand{\pageformat}[6]{\setlength{\hoffset}{-1in}
                  \setlength{\voffset}{-1in}
                  \addtolength{\hoffset}{#5}
                            \addtolength{\voffset}{#6}
                            \setlength{\oddsidemargin}{#1}
                            \setlength{\evensidemargin}{#2}
                            \setlength{\textwidth}{\paperwidth}
                  \addtolength{\textwidth}{-\oddsidemargin}
                  \addtolength{\textwidth}{-\evensidemargin}
                  \addtolength{\textwidth}{-\marginparsep}
                  \addtolength{\textwidth}{-\marginparwidth}
                            \setlength{\topmargin}{#3}
                            \setlength{\textheight}{\paperheight}
                  \addtolength{\textheight}{-\topmargin}
                  \addtolength{\textheight}{-\headheight}
                  \addtolength{\textheight}{-\headsep}
                  \addtolength{\textheight}{-\footskip}
                  \addtolength{\textheight}{-#4}}
\pageformat{2cm}{3cm}{25mm}{25mm}{1pt}{0pt}

\usepackage{ifthen}
\newboolean{article}
    \setboolean{article}{true}
\newboolean{report}
\newboolean{book}
\newboolean{letter}
\newboolean{german}
\newboolean{italian}
\newboolean{nobaselinestretch}
\newboolean{nosectionappendix}
\newboolean{oldtoc}
\newboolean{nosectionequn}
\newboolean{notheorem}

\ifthenelse{\boolean{german}}{
    \usepackage{german}}{}

\usepackage[latin1]{inputenc}

\usepackage{amsmath}
\usepackage{amssymb}
\usepackage[mathscr]{eucal}

\ifthenelse{\boolean{notheorem}}{}{
    \usepackage{theorem}}

\ifthenelse{\boolean{nobaselinestretch}}{}{
    \renewcommand{\baselinestretch}{1.25}}

\newenvironment{env}[2]{\begin{#1}#2\end{#1}}{}
    \newcommand{\beq}[1]{\begin{env}{equation}{#1}}
    \newcommand{\beqn}[1]{\begin{env}{equation*}{#1}}
    \newcommand{\bal}[1]{\begin{env}{align}{#1}}
    \newcommand{\baln}[1]{\begin{env}{align*}{#1}}
    \newcommand{\bga}[1]{\begin{env}{gather}{#1}}
    \newcommand{\bgan}[1]{\begin{env}{gather*}{#1}}
    \newcommand{\bflal}[1]{\begin{env}{flalign}{#1}}
    \newcommand{\bflaln}[1]{\begin{env}{flalign*}{#1}}
    \newcommand{\bmu}[1]{\begin{env}{multline}{#1}}
    \newcommand{\bmun}[1]{\begin{env}{multline*}{#1}}
    \newcommand{\bsp}[1]{\begin{env}{split}{#1}}

    \newcommand{\eeq}{\end{env}}
    \newcommand{\eeqn}{\end{env}}
    \newcommand{\eal}{\end{env}}
    \newcommand{\ealn}{\end{env}}
    \newcommand{\ega}{\end{env}}
    \newcommand{\egan}{\end{env}}
    \newcommand{\eflal}{\end{env}}
    \newcommand{\eflaln}{\end{env}}
    \newcommand{\emu}{\end{env}}
    \newcommand{\emun}{\end{env}}
    \newcommand{\esp}{\end{env}}

\newcommand{\lf}{\vspace{2ex}}

\newcommand{\SQ}[1]{\begin{quotation}\renewcommand{\baselinestretch}{1}\small\noindent\ignorespaces#1\end{quotation}}
\renewcommand{\bf}[1]{\textbf{#1}}
\renewcommand{\it}[1]{\textit{#1}}

\renewcommand{\sf}[1]{\textsf{#1}}

\renewcommand{\tt}[1]{\texttt{#1}}
\newcommand{\hl}[1]{\bf{\it{#1}}}

\newcommand{\mbf}[1]{\mathbf{#1}}
\newcommand{\msf}[1]{\text{\small$\sf{#1}$}}

\newcommand{\cmc}[1]{\mathcal{#1}}
\newcommand{\eus}[1]{\mathscr{#1}}
\newcommand{\euf}[1]{\mathfrak{#1}}
\newcommand{\bb}[1]{\mathbb{#1}}

\newcommand{\mscriptsize}[1]{{\setlength{\arraycolsep}{.3ex}\text{\scriptsize$#1$}}}
\newcommand{\mtiny}[1]{{\setlength{\arraycolsep}{.3ex}\text{\tiny$#1$}}}
\newcommand{\nbd}[1]{$#1$\nobreakdash--}
\newcommand{\ol}[1]{\overline{#1}}

\newcommand{\vt}{\vartheta}

\newcommand{\vp}{\varphi}
\newcommand{\om}{\omega}

\newcommand{\abs}[1]{\left\lvert#1\right\rvert}
\newcommand{\norm}[1]{\left\lVert#1\right\rVert}

\newcommand{\bfam}[1]{\bigl(#1\bigr)}

\newcommand{\AB}[1]{\langle#1\rangle}
\newcommand{\bAB}[1]{\bigl\langle#1\bigr\rangle}
\newcommand{\BAB}[1]{\Bigl\langle#1\Bigr\rangle}
\newcommand{\CB}[1]{\{#1\}}
\newcommand{\bCB}[1]{\bigl\{#1\bigr\}}
\newcommand{\BCB}[1]{\Bigl\{#1\Bigr\}}
\newcommand{\SB}[1]{[#1]}

\newcommand{\Matrix}[1]{\begin{pmatrix}#1\end{pmatrix}}

\newcommand{\sMatrix}[1]{\mscriptsize{\Matrix{#1}}}
\newcommand{\tMatrix}[1]{\mtiny{\Matrix{#1}}}

\newcommand{\set}[2][]{
    \ifthenelse{\equal{#1}{}}{
        \CB{#2}}{
        \CB{#1~|~#2}}}
\newcommand{\bset}[2][]{
    \ifthenelse{\equal{#1}{}}{
        \bCB{#2}}{
        \bCB{#1~|~#2}}}
\newcommand{\Bset}[2][]{
    \ifthenelse{\equal{#1}{}}{
        \BCB{#2}}{
        \BCB{#1~\big|~#2}}}

\newlength{\thmboxwidth}
\newcommand{\thmbox}[2][14cm]{
\setlength{\thmboxwidth}{#1}
\begin{center}\fbox{\begin{minipage}{\thmboxwidth}
#2
\end{minipage}}\end{center}
}

\DeclareMathOperator{\ls}{\normalfont\msf{span}}
\DeclareMathOperator{\cls}{\ol{\ls}}

\DeclareMathOperator{\id}{\normalfont\msf{id}}

\newcommand{\C}{\bb{C}}

\newcommand{\E}{\bb{E}}

\newcommand{\N}{\bb{N}}

\newcommand{\R}{\bb{R}}

\newcommand{\cA}{\cmc{A}}
\newcommand{\cB}{\cmc{B}}
\newcommand{\cC}{\cmc{C}}
\newcommand{\cD}{\cmc{D}}

\newcommand{\cS}{\cmc{S}}

\newcommand{\sB}{\eus{B}}

\newcommand{\sE}{\eus{E}}
\newcommand{\sF}{\eus{F}}

\newcommand{\sK}{\eus{K}}
\newcommand{\sL}{\eus{L}}

\newcommand{\sN}{\eus{N}}

\newcommand{\sS}{\eus{S}}

\newcommand{\ek}{\euf{k}}
\newcommand{\el}{\euf{l}}

\newcommand{\et}{\euf{t}}

\newcommand{\eK}{\euf{K}}
\newcommand{\eL}{\euf{L}}

\newcommand{\eT}{\euf{T}}

\newcommand{\U}{\mbf{1}}

\newcommand{\G}{\Gamma}

\newcommand{\I}{{I\!\!\!\;I}}

\newcommand{\s}{\text{\scriptsize$\sS$}}

\ifthenelse{\boolean{nosectionequn}}{}{
    \numberwithin{equation}{section}
    }

\ifthenelse{\boolean{article}\or\boolean{letter}\or\boolean{nosectionequn}}{
    \setboolean{nosectionappendix}{true}}{}
\ifthenelse{\boolean{nosectionappendix}}{}{
    \renewcommand{\appendix}{
        \chapter*{\appendixname}
        \addcontentsline{toc}{chapter}{\appendixname}
        \renewcommand{\thesection}{\Alph{section}}
        \setcounter{section}{0}}}
   
\ifthenelse{\boolean{report}\or\boolean{book}}{
    }{}

\ifthenelse{\boolean{notheorem}}{}{
        \newcommand{\notename}{Note.}
        \newcommand{\mnname}{Mathematical note.}
        \newcommand{\enname}{End of the note.}
        \newcommand{\definame}{Definition.}
        \newcommand{\propname}{Proposition.}
        \newcommand{\lemname}{Lemma.}
        \newcommand{\exname}{Example.}
        \newcommand{\exername}{Exercise.}
        \newcommand{\remname}{Remark.}
        \newcommand{\obname}{Observation.}
        \newcommand{\thmname}{Theorem.}
        \newcommand{\corname}{Corollary.}
        \newcommand{\proofname}{Proof.}
    \ifthenelse{\boolean{german}}{
        \renewcommand{\mnname}{Mathematische Notiz.}
        \renewcommand{\enname}{Ende der Notiz.}
        \renewcommand{\exname}{Beispiel.}
        \renewcommand{\exername}{Übung.}
        \renewcommand{\remname}{Bemerkung.}
        \renewcommand{\obname}{Beobachtung.}
        \renewcommand{\thmname}{Satz.}
        \renewcommand{\corname}{Korollar.}
        \renewcommand{\proofname}{Beweis.}}{}
    \ifthenelse{\boolean{italian}}{
        \renewcommand{\mnname}{Nota matematica.}
        \renewcommand{\enname}{Fina della nota.}
        \renewcommand{\definame}{Definizione.}
        \renewcommand{\propname}{Proposizione.}
        \renewcommand{\exname}{Esempio.}
        \renewcommand{\exername}{Esercizio.}
        \renewcommand{\remname}{Nota.}
        \renewcommand{\obname}{Osservazione.}
        \renewcommand{\thmname}{Teorema.}
        \renewcommand{\corname}{Corollario.}
        \renewcommand{\proofname}{Dimostrazione.}

       \renewcommand{\appendixname}{Appendice}

       }{}
    \theoremheaderfont{\normalfont\bfseries}
    \theoremstyle{change}
        \theorembodyfont{\rmfamily}
            \newtheorem{emp}{}[section]
                \newcommand{\bemp}[1][]{
                    \begin{emp}\hskip-\labelsep\bf{#1}\hskip\labelsep}
                \newcommand{\eemp}{\end{emp}}
\newtheorem{itemp}[emp]{}
                \newcommand{\bitemp}[1][]{
                    \begin{itemp}\hskip-\labelsep\bf{#1}\hskip\labelsep\normalfont\itshape}
                \newcommand{\eitemp}{\end{itemp}}
            \newtheorem{note}[emp]{\notename}
                \newcommand{\bnote}{\begin{note}}
                \newcommand{\enote}{\end{note}}
            \newtheorem{mn}[emp]{\mnname}
                \newcommand{\bnm}{\begin{mn}~\begin{quotation}\renewcommand{\baselinestretch}{1}\small\noindent\ignorespaces}
                \newcommand{\enm}{\end{quotation}\hfill\bf{\enname}\end{mn}}
            \newtheorem{ex}[emp]{\exname}
                \newcommand{\bex}{\begin{ex}}
                \newcommand{\eex}{\end{ex}}
            \newtheorem{exer}[emp]{\exername}
                \newcommand{\bexer}{\begin{exer}}
                \newcommand{\eexer}{\end{exer}}
            \newtheorem{defi}[emp]{\definame}
                \newcommand{\bdefi}{\begin{defi}}
                \newcommand{\edefi}{\end{defi}}
            \newtheorem{rem}[emp]{\remname}
                \newcommand{\brem}{\begin{rem}}
                \newcommand{\erem}{\end{rem}}
            \newtheorem{ob}[emp]{\obname}
                \newcommand{\bob}{\begin{ob}}
                \newcommand{\eob}{\end{ob}}
        \theorembodyfont{\normalfont\itshape}
            \newtheorem{thm}[emp]{\thmname}
                \newcommand{\bthm}{\begin{thm}}
                \newcommand{\ethm}{\end{thm}}
            \newtheorem{prop}[emp]{\propname}
                \newcommand{\bprop}{\begin{prop}}
                \newcommand{\eprop}{\end{prop}}
            \newtheorem{cor}[emp]{\corname}
                \newcommand{\bcor}{\begin{cor}}
                \newcommand{\ecor}{\end{cor}}
            \newtheorem{lem}[emp]{\lemname}
                \newcommand{\blem}{\begin{lem}}
                \newcommand{\elem}{\end{lem}}
\newenvironment{empn}[1]{\lf\noindent\bf{#1}\ignorespaces\hskip\labelsep}{\lf}
		\newcommand{\bempn}[1]{\begin{empn}{#1}}
		\newcommand{\eempn}{\end{empn}}
		\newcommand{\bitempn}[1]{\begin{empn}{#1}\normalfont\itshape}
		\newcommand{\eitempn}{\end{empn}}
                \newcommand{\bnmn}{\begin{empn}{\mnname}~\begin{quotation}\renewcommand{\baselinestretch}{1}\small\noindent\ignorespaces}
                \newcommand{\enmn}{\end{quotation}\hfill\bf{\enname}\end{empn}}
		\newcommand{\bexn}{\begin{empn}{\exname}}
		\newcommand{\eexn}{\end{empn}}
		\newcommand{\bexern}{\begin{empn}{\exername}}
		\newcommand{\eexern}{\end{empn}}
		\newcommand{\bdefin}{\begin{empn}{\definame}}
		\newcommand{\edefin}{\end{empn}}
		\newcommand{\bremn}{\begin{empn}{\remname}}
		\newcommand{\eremn}{\end{empn}}
		\newcommand{\bobn}{\begin{empn}{\obname}}
		\newcommand{\eobn}{\end{empn}}

\newcommand{\qedsymbol}{~\rule[-0.35mm]{2mm}{2mm}}
    \newcounter{proof}[emp]
    \newenvironment{Proof}[1]{
        \vspace{1ex}
        \renewcommand{\item}[1][\stepcounter{proof}(\roman{proof})]%
            {##1\hskip\labelsep}
        \noindent\textsc{#1\hskip\labelsep}}{
        \nolinebreak\qedsymbol}
    \newcommand{\proof}[1][\proofname]{
        \begin{Proof}{#1}\ignorespaces}
    \newcommand{\qed}{\end{Proof}}
    \newcommand{\noqed}{
        \renewcommand{\qedsymbol}{}
        \end{Proof}}}
    \ifthenelse{\boolean{italian}}{
        \renewcommand{\proofname}{Dimostrazione.}}{}

\usepackage[varg]{txfonts}

\usepackage[hypertex]{hyperref}

\begin{document}

\bibliographystyle{amsalpha}

\title{Hilbert Modules---Square Roots of Positive Maps\thanks{This work has been supported by research funds of the Dipartimento S.E.G.e S.\ of University of Molise and of the Italian MUR (PRIN 2007).}}

\author{Michael Skeide}

\date{June 2009}

\maketitle

\begin{abstract}
\noindent
We reflect on the notions of positivity and square roots. We review many examples which underline our thesis that square roots of positive maps related to \nbd{*}algebras are Hilbert modules. As a result of our considerations we discuss requirements a notion of positivity on a \nbd{*}algebra should fulfill and derive some basic consequences.
\end{abstract}

\section{Introduction}\label{intro}

Let $S$ denote a set, and let $\ek$ denote a map $S\times S\rightarrow\C$. Everybody knows that such a \hl{kernel} $\ek$ \hl{over} $S$ is called \hl{positive definite} if
\beq{\label{PDc}
\sum_{i,j}\bar{z}_i\ek^{\sigma_i,\sigma_j}z_j
~\ge~
0
}\eeq
for all finite choices of $\sigma_i\in S$ and $z_i\in\C$.

What is is the best way to show that some thing $x$ is positive? The best way is writing $x$ as a \hl{square}! It would, then, be justified to call an object $y$ the positive thing $x$'s \hl{square root}, if by writing down the object $y$'s square we get back $x$. By \it{square}, of course, we mean a \hl{complex square} like $\bar{y}y$ ($y$ a complex number) or $y^*y$ ($y$ in a \nbd{C^*}algebra) or $\AB{y,y}$ ($y$ being in a Hilbert space).

Of course, for each choice $\sigma_i\in S$ and $z_i\in\C$ we may calculate the positive number in \eqref{PDc} and write down its positive square root $p(\sigma_1,\ldots,\sigma_n,z_1,\ldots,z_n)$ (or any other complex square root) and that's it. Although, the collection of all $p$ contains the full information about $\ek$ (for instance by suitable polarization procedures or by differentiation with respect to the parameters $z_i$), it is uncomfortable to do that. Also, the knowledge of some $p(\sigma_1,\ldots,\sigma_n,z_1,\ldots,z_n)$ for a fixed choice, does not at all help computing $p(\sigma_1,\ldots,\sigma_{n-1},z_1,\ldots,z_{n-1})$ for the same choice. We gain a bit but not very much, if we calculate for each choice $\sigma_1,\ldots,\sigma_n$ the positive (or some other) square root $P(\sigma_1,\ldots,\sigma_n)\in M_n$ of the positive matrix $\bfam{\ek^{\sigma_i\sigma_j}}_{i,j}\in M_n$. Still, the knowledge of some $P(\sigma_1,\ldots,\sigma_n)$ for a certain choice does not help computing $P(\sigma_1,\ldots,\sigma_{n-1})$ for the same choice. (Exercise: Try it and explain why it does not help!)

We wish something that allows easily to recover the function $\ek$ and that still gives evidence of positivity of the expressions in \eqref{PDc} by writing them as square. The solution to that problem is the well-known \hl{Kolmogorov decomposition}.

\bthm\label{KDcthm}
For every positive definite kernel $\ek$ over a set $S$ with values in $\C$ there exist a Hilbert space $H$ and a map $i\colon S\rightarrow H$ such that
\beqn{
\ek^{\sigma,\sigma'}
~=~
\AB{i(\sigma),i(\sigma')}
}\eeqn
for all $\sigma,\sigma'\in S$.
\ethm

\proof
On the vector space $S_\C:=\bigoplus_{\sigma\in S}\C=\Bset[~\bfam{z_\sigma}_{\sigma\in S}]{\#\CB{\sigma\colon z_\sigma\ne0}<\infty~}$ we define a sesquilinear form
\beqn{
\BAB{\bfam{z_\sigma}_{\sigma\in S},\bfam{z'_\sigma}_{\sigma\in S}}
~:=~
\sum_{\sigma,\sigma'\in S}\bar{z}_\sigma\ek^{\sigma,\sigma'}z'_{\sigma'}.
}\eeqn
Since $\ek$ is positive definite, this form is positive. Denote $e_\sigma:=\bfam{\delta_{\sigma,\sigma'}}_{\sigma'\in S}$. Then $\AB{e_\sigma,e_{\sigma'}}=\ek^{\sigma,\sigma'}$. Denote by $H$ the \hl{Hausdorff completion} of $\bigoplus_{\sigma\in S}\C$ (that is, quotient out the subspace $\sN$ of length-zero elements and complete that pre-Hilbert space). Then $H$ with the function $i$ defined by $i(\sigma):=e_\sigma+\sN$ has the claimed properties.\qed

\lf
Note that the subset $i(S)$ of $H$ as constructed in the proof is total. Therefore, the pair $(H,i)$ has the following universal property: If $(G,j)$ is another Kolmogorov decomposition of $\ek$, then there is a unique bounded linear operators $v\colon H\rightarrow G$ such that $vi(\sigma)=j(\sigma)$ for all $\sigma\in S$. Note that $v$ is isometric so that $(H,i)$ is determined by that universal property up to unique unitary equivalence. This is just the same as the square root $p$ of a positive number $k$, which is determined up to a unitary operator $e^{i\vp}$ on the one-dimensional Hilbert space $\C$.

We like to think of the \hl{minimal Kolmogorov} construction $(H,i)$ as \bf{the} square root of the kernel $\ek$. Obviously, every Hilbert space arises in that way. (Simply take the kernel $\ek^{h,h'}:=\AB{h,h'}$. Then $(H,i\colon h\mapsto h)$ has the universal property.)

It is the scope of these notes to establish the idea of inner product spaces (like Hilbert modules) as square roots of maps that are positive in some sense. Apart from many instances of this interpretation, we intend also to discuss the just mentioned uniqueness issue for square roots, and to present the rudiments of what we consider a ``good'' notion of positivity in \nbd{*}algebras:
\thmbox{
In a ``good notion of positivity'' it should be a theorem that all positive things have a sort of square root.
}
Another scope is to point out the following insight about composition of positive things.
\thmbox{
In the noncommutative world, if one wishes to compose positive things to get new ones, then these positive things must be maps on \nbd{*}algebras, not elements in \nbd{*}algebras.
}
Many of our examples have to do with product systems. We should mention that we systematically omit mentioning any relationship that has to do with commutants of von Neumann correspondences (Skeide \cite{Ske03c}). We refer the interested reader to the survey Skeide \cite{Ske08a}.

\section{Kernels with values in a $C^*$--algebra}\label{PDSEC}

If $\ek\colon S\times S\rightarrow\cB$ is a kernel \hl{over} $S$ with values in a \nbd{C^*}algebra $\cB$, then everything  goes precisely as in the scalar-valued case, just that now the space emerging by Kolmogorov decomposition is a Hilbert \nbd{\cB}module.

A kernel $\ek$ is \hl{positive definite} (or a \hl{PD-kernel}) if
\beq{\label{PD}
\sum_{i,j}b_i^*\ek^{\sigma_i,\sigma_j}b_j
~\ge~
0
}\eeq
for all finite choices of $\sigma_i\in S$ and $b_i\in\C$. Let us equip the right \nbd{\cB}module $E_0:=S_\C\otimes\cB=\bigoplus_{\sigma\in S}\cB=\Bset[~\bfam{b_\sigma}_{\sigma\in S}]{\#\CB{\sigma\colon b_\sigma\ne0}<\infty~}$ with the sesquilinear map $\AB{\bullet,\bullet}\colon E_0\times E_0\rightarrow\cB$
\beqn{
\BAB{\bfam{b_\sigma}_{\sigma\in S},\bfam{b'_\sigma}_{\sigma\in S}}
~:=~
\sum_{\sigma\in S}b_\sigma^*\ek^{\sigma,\sigma'}b'_{\sigma'}.
}\eeqn
Equation \eqref{PD} is born to to make $\AB{\bullet,\bullet}$ \hl{positive}: $\AB{x,x}\ge0$ for all $x\in E_0$. It also is \hl{right \nbd{\cB}linear}: $\AB{x,yb}=\AB{x,y}b$ for all $x,y\in E_0$ and $b\in\cB$. In other words, $\AB{\bullet,\bullet}$ is a \hl{semiinner product} and $E_0$ is a \hl{semi-Hilbert \nbd{\cB}module}. By making appropriate use of \hl{Cauchy-Schwarz inequality}
\beqn{
\AB{x,y}\AB{y,x}
~\le~
\norm{\AB{y,y}}\AB{x,x}
}\eeqn
(Paschke \cite{Pas73}), the function $x\mapsto\sqrt{\norm{\AB{x,x}}}$ is a seminorm. So, we may divide out the right submodule of length-zero elements $\sN$. In other words, $E_0/\sN$ is a \hl{pre-Hilbert \nbd{\cB}module}, that is, it is a  semi-Hilbert \nbd{\cB}module where $\AB{x,x}=0$ implies $x=0$ for all $x\in E_0/\sN$. Moreover, $\norm{xb}\le\norm{x}\norm{b}$ so that we may complete the quotient. In other words,  $E:=\ol{E_0/\sN}$ is a \hl{Hilbert \nbd{\cB}module}, that is, $E$ is a complete pre-Hilbert \nbd{\cB}module.

Recall that $e_\sigma\otimes b=\bfam{\delta_{\sigma,\sigma'}b}_{\sigma'\in S}$. If $\cB$ is unital, then $i(\sigma):=e_\sigma\otimes\U+\sN$ fulfills $\AB{i(\sigma),i(\sigma')}=\ek^{\sigma,\sigma'}$ and $\cls i(S)\cB=E$. If $\cB$ is nonunital, then choose an approximate unit $\bfam{u_\lambda}_{\lambda\in\Lambda}$ for $\cB$, and verify that $\bfam{e_\sigma\otimes u_\lambda+\sN}_{\lambda\in\Lambda}$ is a Cauchy net in $E$. Define $i(\sigma):=\lim_\lambda e_\sigma\otimes u_\lambda+\sN$. In conclusion:

\bthm\label{KdBthm}
If $\ek$ is a \nbd{\cB}valued PD-kernel over $S$, then there is a pair $(E,i)$ of a Hilbert \nbd{\cB}module $E$ and map $i\colon S\rightarrow E$ satisfying
\beqn{
\AB{i(\sigma),i(\sigma')}
~=~
\ek^{\sigma,\sigma'}
}\eeqn
for all $\sigma,\sigma'\in S$ and $\cls i(S)\cB=E$. Moreover, if $(F,j)$ is another pair fulfilling $\AB{j(\sigma),j(\sigma')}=\ek^{\sigma,\sigma'}$, then the map $i(\sigma)\mapsto j(\sigma)$ extends to a unique \hl{isometry} (that is, an inner product preserving map) $E\rightarrow F$.
\ethm

By the universal property, it follows that the pair $(H,i)$ is determined up to unique unitary equivalence. (A \hl{unitary} is a surjective isometry.) We refer to it as the \hl{minimal Kolmogorov decomposition} of $\ek$.

Once more, every Hilbert module $E$ arises in that way, as the Kolmogorov decomposition $(E,\id_E)$ of the PD-kernel $(x,y)\mapsto\AB{x,y}$ over $E$. We, therefore, like to think of Hilbert modules as square roots of PD-kernels.

\bex\label{1ptex}
For a positive element $b\in\cB$ we may define the PD-kernel $\ek\colon(\om,\om)\mapsto b$ over the one-point set $S=\CB{\om}$. If we choose an element $\beta\in\cB$ such that $\beta^*\beta=b$, then the right ideal $E=\ol{\beta\cB}$ generated by $\beta$ with inner product $\AB{x,y}:=x^*y$ is a Hilbert \nbd{\cB}module. Moreover, the map $i\colon\om\mapsto\beta$ fulfills $\AB{i(\om),i(\om)}=\ek^{\om,\om}$ and $\cls i(\om)\cB=E$.

If $\beta'$ is a another square root, then Theorem \ref{KdBthm} tells us that $\beta\mapsto\beta'$ extends as a unitary from $E=\ol{\beta\cB}$ to $E'=\ol{\beta'\cB}$. But more cannot be said about different choices of square roots. For instance, if $b=\U$, then every isometry $v\in\cB$ is a square root. But as subsets of $\cB$ the sets $v\cB(=\ol{v\cB})$ can be quite different. It can be all $\cB$. (This happens if and only if $v$ is a unitary.) But if $v$ and $v'$ fulfill $v^*v'=0$, then they are even orthogonal to each other. If $\cB$ is unital and $\ol{\beta\cB}=\cB$, then $\beta$ is necessarily invertible. (Exercise!) If $\beta'\in\cB$ fulfills $\beta^*\beta=\beta'^*\beta'$, then it is easy to show that $\beta'\beta^{-1}$ is a unitary. $(\beta^{-1})^*\beta'^*\beta'\beta^{-1}=(\beta^*)^{-1}\beta^*\beta\beta^{-1}=\U$.

Only, the picture of Kolmogorov decomposition for the kernel on a one-point set allows to make a precise statement.
\eex

\SQ{\vspace{-3ex}
\bnote
It seems that the concept of PD-kernels with values in an abstract \nbd{C^*}algebra has not been considered before Barreto, Bhat, Liebscher, and Skeide \cite{BBLS04}. The classical Stinespring theorem \cite{Sti55} for CP-maps with values in a concrete \nbd{C^*}algebra $\cB\subset\sB(G)$ is proved by using the Kolmogorov decomposition for a \nbd{\C}valued PD-kernel; see Remark \ref{Stinerem}. However, analogue constructions for CP-maps with values in $\sB^a(F)$ (the algebra of adjointable operators on a Hilbert \nbd{\cC}modules $F$) by Kasparov \cite{Kas80} and Lance \cite{Lan95} use proofs similar to Paschke's GNS-construction \cite{Pas73} for CP-maps; see Note \ref{KSGNSnote}. Closest is Murphy's result in\cite{Mur97} for \nbd{\sB^a(F)}valued kernels, whose proof uses techniques like reproducing kernels (Aronszajin \cite{Aro50}); see Szafraniec's survey \cite{Sza09p}.
\enote
}

\section{Composing PD-kernels?}\label{PDcompSEC}

It is well known that two positive definite \nbd{\C}valued kernels $\el$ and $\ek$ over the same set $S$ may composed by \hl{Schur product}, that is, by the pointwise product
\beqn{
(\el\ek)^{\sigma,\sigma'}
~:=~
\el^{\sigma,\sigma'}\ek^{\sigma,\sigma'},
}\eeqn
and the result is again a PD-kernel over $S$. Note that this Schur product of \nbd{\C}valued kernels is commutative.

Of course, we may define the Schur product of \nbd{\cB}valued kernels by the same formula . But now the product, in general, depends on the order. However, if $\el\ek\ne\ek\el$, then neither of the two products  is PD. (Note that by Kolmogorov decomposition, a PD-kernel is necessarily \hl{hermitian}: ${\ek^{\sigma,\sigma'}}^*=\AB{i(\sigma),i(\sigma')}^*=\AB{i(\sigma'),i(\sigma)}=\ek^{\sigma',\sigma}$.)

Does it help if we try to compose the square roots? Let us choose two positive elements $b=\beta^*\beta$ and $c=\gamma^*\gamma$, and, as in Example \ref{1ptex}, consider the two PD-kernels $\ek\colon(\om,\om)\mapsto b$ and $\el\colon(\om,\om)\mapsto c$ over the one-point set $S=\CB{\om}$. We may take the two square roots $\beta$ and $\gamma$, multiply them, and use their product $\beta\gamma$ to define a PD-kernel $(\om,\om)\mapsto(\beta\gamma)^*(\beta\gamma)=\gamma^*\beta^*\beta\gamma$ on $S$.

There are two things to be noted. First, if $\beta$ and $\gamma$ do note commute, then the ``composed'' kernel depends on the order. This as such is not too disturbing in a noncommutative context.

\SQ{\vspace{-3ex}
\bnote
Bercovici \cite{Ber05} and Franz \cite{Fra09} use such a procedure of a product in the definition of multiplicative monotone convolution of probability measures.
\enote
}
Second, and much more crucial, the kernel $\el$ alone does not allow to determine that ``composition''. Or the other way round, different square roots $\gamma$ of $\el$ do, in general, not give rise to the same composition. What we know about the kernel is equivalently coded in its minimal Kolmogorov decomposition. However, as pointed out in Example \ref{1ptex}, different square roots $\gamma$ are indistinguishable both from the point of view of Kolmogorov decomposition and from the point of view of the kernel itself.

The puzzle is resolved, if we observe that ,actually, we have to compute the map $\gamma^*\bullet\gamma\colon b\mapsto\gamma^*b\gamma$ --- a map with strong positivity properties. If we wish to compose $\el$ with arbitrary kernels $\ek$, then we need the entire information about that map. That information is encoded in the left ideal generated by $\gamma$. Doing also here a Kolmogorov type construction, we end up with the two-sided ideal, that is, the Hilbert \nbd{\cB}bimodule, generated by $\gamma$; see Example \ref{1pKex} below.

What we just discussed for a one-point set, for general sets $S$ gives rise to the notion of completely positive definite (CPD) kernels. The Kolmogorov decomposition for CPD-kernels will result in a Hilbert bimodule rather than in an Hilbert module.
\thmbox{
CPD-kernels may be composed, and the Kolmogorov decomposition for the composition of two CPD-kernels is reflected by the Kolmogorov decompositions of the factors. CPD-kernels are, therefore, the ``correct'' generalization of \nbd{\C}valued PD-kernels.
}
This will be subject of the next section.

\section{CPD-kernels}\label{CPDSEC}

In a noncommutative context we have seen that, if we wish to compose kernels fulfilling some positivity condition in a way that preserves positivity, then it is practically forced to switch from kernels with values in $\cB$ to kernels with values in the bounded maps on $\cB$. Once we have map-valued kernels, there is no longer a reason that domain and codomain must coincide.

\bdefi
Let $S$ be set and let $\eK\colon S\times S\rightarrow\sB(\cA,\cB)$ be a kernel over $S$ with values in the bounded maps from a \nbd{C^*}algebra $\cA$ to a \nbd{C^*}algebra $\cB$. We say $\eK$ is a \hl{completely positive definite kernel} (or \hl{CPD-kernel}) over $S$ from $\cA$ to $\cB$ if
\beq{\label{CPD}
\sum_{i,j}b_i^*\eK^{\sigma_i,\sigma_j}(a_i^*a_j)b_j
~\ge~
0
}\eeq
for all finite choices of $\sigma_i\in S,a_i\in\cA,b_i\in\cB$. If $\cA=\cB$, then we also say a kernel \hl{on} $\cB$.
\edefi

There are two possibilities to find the appropriate Kolmogorov decomposition for CPD-kernels. As pointed out in the end of the last section, it is no surprise that we obtain a Hilbert bimodule or, more fashionably, a \it{correspondence}. Recall that a \hl{correspondence} from $\cA$ to $\cB$ is a Hilbert \nbd{\cB}module $E$ with a nondegenerate(\bf{!}) left action of $\cA$ such that $\AB{x,ay}=\AB{a^*x,y}$ for all $x,y\in E$ and $a\in\cA$.

\bthm\label{KdABthm}
If $\eK$ is a CPD-kernel over $S$ from a unital \nbd{C^*}algebra $\cA$ to a \nbd{C^*}algebra $\cB$, there is pair $(E,i)$ consisting of a correspondence $E$ from $\cA$ to $\cB$ and a map $i\colon S\rightarrow E$ such that
\beqn{
\AB{i(\sigma),ai(\sigma')}
~=~
\eK^{\sigma,\sigma'}(a)
}\eeqn
for all $\sigma,\sigma'\in S$ and $a\in\cA$, and such that $E=\cls\cA i(S)\cB$. Moreover, if $(F,j)$ is another pair fulfilling $\AB{j(\sigma),aj(\sigma')}=\eK^{\sigma,\sigma'}(a)$, then the map $i(\sigma)\mapsto j(\sigma)$ extends to a unique bilinear isometry $E\rightarrow F$.
\ethm
 
We refer to $(E,i)$ as the \hl{Kolmogorov decomposition} of the CPD-kernel $\eK$. By the universal property stated in the theorem, it is uniquely determined up to bilinear unitary equivalence. We also shall refer to $E$ as the \hl{GNS-correspondence} and to $i$ as the \hl{cyclic map}; see Note \ref{GNSnote}.

\proof[Proof of Theorem \ref{KdABthm}.~~]
First possibility: By \eqref{CPD} it immediately follows that the kernel $\ek^{(a,\sigma),(a',\sigma')}:=\eK^{\sigma,\sigma'}(a^*a')$ over $\cA\times S$ is positive definite. Denote by $(E,\tilde{i})$ its Kolmogorov decomposition according to Theorem \ref{KdBthm}. On the subset $\tilde{i}(\cA\times S)$ define a left action by $a\tilde{i}(a',\sigma'):=\tilde{i}(aa',\sigma')$. This action fulfills $\AB{\tilde{i}(a',\sigma'),a\tilde{i}(a'',\sigma'')}=\AB{a^*\tilde{i}(a',\sigma'),\tilde{i}(a'',\sigma'')}$. Recall that the set $\tilde{i}(\cA\times S)$ generates $E$ as Hilbert \nbd{\cB}module. It is easy to prove that an action fulfilling the \nbd{*}condition on a generating subset of $E$ extends well-defined and uniquely to a left action on all of $E$. The pair $(E,i)$ satisfies the stated properties.

Second possibility: Instead of appealing to Theorem \ref{KdBthm}, we imitate its proof. Indeed, if we equip the \nbd{\cA}\nbd{\cB}bimodule $E_0:=\cA\otimes S_\C\otimes\cB$ with the sesquilinear map defined by setting
\beqn{
\AB{a\otimes e_\sigma\otimes b,a'\otimes e_{\sigma'}\otimes b'}
~:=~
b^*\eK^{\sigma,\sigma'}(a^*a')b',
}\eeqn
then the condition in Equation \eqref{CPD} is born to make it an semiinner product, which also fulfills $\AB{a'\otimes e_{\sigma'}\otimes b',aa''\otimes e_{\sigma''}\otimes b''}=\AB{a^*a'\otimes e_{\sigma'}\otimes b',a''\otimes e_{\sigma''}\otimes b''}$. We divide out the \nbd{\cA}\nbd{\cB}submodule of length-zero elements and complete, obtaining that way an \nbd{\cA}\nbd{\cB}correspondence $E$. The map $i(\sigma):=\lim_\lambda\U_\cA\otimes e_\sigma\otimes u_\lambda+\sN$ completes the construction.\qed

\SQ{\vspace{-3ex}
\bnote\label{GNSnote}
For one-point sets $S=\CB{\om}$ we get back the definition of CP-maps between \nbd{C^*}algebras, and the second proof of Theorem \ref{KdABthm} is just Paschke's \hl{GNS-construction} for CP-maps; see \cite{Pas73}.
\enote
}

\brem\label{Stinerem}
Why did we present two proofs for Theorem \ref{KdABthm}? The first proof is more along classical lines: From the input data write down some kernel (in classical applications almost always \nbd{\C}valued), show it is positive definite, and do the Kolmogorov decomposition. Only then start working in order to show that this Hilbert module (or, usually, Hilbert space in classical applications) has the desired structure.

For instance, classical proofs of the Stinespring construction for a CP-map $T\colon\cA\rightarrow\cB$ work approximately like that (cf.\ also Example \ref{Stinespring}): Represent your \nbd{C^*}algebra $\cB$ faithfully on a Hilbert space $G$ and define a \nbd{\C}valued kernel $\ek$ over the set $\cA\times G$ as $\ek^{(a,g),(a',g')}:=\AB{g,T(a^*a')g'}$. Work in order to prove it is positive definite. Do the Kolmogorov decomposition to get the pair $(H,\tilde{i})$. Work again in order to show that $\rho(a)\colon\tilde{i}(a',g)\mapsto\tilde{i}(aa',g)$ determines a representation $\rho$ of $\cA$ on $H$. Work still more in order to show that $v\colon g\mapsto i(\U,g)$ defines a bounded map such that $v^*\rho(a)v=T(a)$ for all $a\in\cA$. But how much work, is'nt it!

The second proof is different. We want an \nbd{\cA}\nbd{\cB}bimodule? That proof starts by writing down the \nbd{\cA}\nbd{\cB}bimodule $E_0$. It is immediate from the input data how to define a semiinner product that turns it into a semicorrespondence from $\cA$ to $\cB$. Apply the generalities from the theory of correspondences (essentially, Cauchy-Schwarz inequality; cf.\ Remark \ref{CSIrem}) that tells you that one may divide out kernels of semiinner products and complete. Identify the elements $i(\sigma)$, which, by definition of the inner product, fulfill the stated property. (Needless to say that representing the \nbd{C^*}algebra $\cB$ is neither necessary nor useful. Anyway, never represent an abstract \nbd{C^*}algebra that it is not given as a concrete operator algebra from the beginning, unless you are going to prove a theorem that is explicitly about representations of $\cB$ on a Hilbert space!)
\erem

After this deviation let us return to our subject: Square roots of positive things. Let us note that the Kolmogorov decomposition is a ``good'' square root of $\eK$. It allows easily to get $\eK$ back as $\eK^{\sigma,\sigma'}(a)=\AB{i(\sigma),ai(\sigma')}$. It puts into immediate evidence why $\eK$ is completely positive definite if it can be recovered by a pair $(E,i)$. Indeed,
\beqn{
\sum_{i,j}b_i^*\eK^{\sigma_i,\sigma_j}(a_i^*a_j)b_j
~=~
\BAB{\sum_ia_ii(\sigma_i)b_i,\sum_ia_ii(\sigma_i)b_i}
~\ge~
0.
}\eeqn
And it is unique up to suitable unitary equivalence of correspondences.

\bex\label{1pKex}
Let us return to the situation with the PD-kernels $\ek\colon(\om,\om)\mapsto b=\beta^*\beta$ and $\el\colon(\om,\om)\mapsto c=\gamma^*\gamma$ over the one-point set $S=\CB{\om}$ as discussed in Section \ref{PDcompSEC}. We have noted that in order to understand the composition defined as $(\om,\om)\mapsto(\beta\gamma)^*(\beta\gamma)=\gamma^*\beta^*\beta\gamma$ whatever $b$ might be, we must know the map $\eL\colon b\mapsto\gamma^*b\gamma$ rather than  just the kernel $\el$. Of course, $\eL$ is a CPD-kernel over $S$, and its Kolmogorov decomposition is $(F,j)$ with $F=\cls\cB\gamma\cB$ and $j(\om)=\gamma$. Theorem \ref{KdABthm} tells us that another $\gamma'$ gives the same $\eL$ if and only if $\gamma\mapsto\gamma'$ extends as an isomorphism of correspondences over $\cB$. This is precisely the case if $(\beta\gamma)^*(\beta\gamma)$ and $(\beta\gamma')^*(\beta\gamma')$ coincide whatever $\beta$ is.

Having replaced  $\el\colon(\om,\om)\mapsto c$ with $\eL\colon(\om,\om)\mapsto\gamma^*\bullet\gamma$, it is only natural to do the same with $\ek\colon(\om,\om)\mapsto b$ and to replace it with $\eK\colon(\om,\om)\mapsto\beta^*\bullet\beta$. We can easily define $\eL\circ\eK\colon(\om,\om)\mapsto(\beta\gamma)^*\bullet(\beta\gamma)$. This brings us directly to the question of composition of kernels that motivated the definition of CPD-kernels.
\eex

\bdefi\label{compdef}
Let $S$ be a set, let $\eK$ be a kernel over $S$ from $\cA$ to $\cB$, and  let $\eL$ be a kernel over $S$ from $\cB$ to $\cC$. Then we define their \hl{composition} or \hl{Schur product} $\eL\circ\eK$ over $S$ from $\cA$ to $\cC$ by pointwise composition, that is, by
\beqn{
(\eL\circ\eK)^{\sigma,\sigma'}
~:=~
\eL^{\sigma,\sigma'}\circ\eK^{\sigma,\sigma'}.
}\eeqn
\edefi

As the proof of the fact that the composition of CPD-kernels is CPD has to do with some of our  considerations about positivity we wish to make in general, we postpone it to the next section. But once this is settled, it is clear what we will understand by a \hl{CPD-semigroup} $\eT=\bfam{\eT_t}_{t\ge0}$ of CPD-kernels $\eT_t$ over a set $S$ from $\cB$ to $\cB$ (or \hl{on} $\cB$).

\SQ{\vspace{-3ex}
\bnote
A single CPD-kernel from $\cA$ to $\cB$ over the finite set $S=\CB{1,\ldots,n}$ has been defined by Heo \cite{Heo99} under the name of \it{completely multi-positive map}. Special CPD-semigroups over $S=\CB{0,1}$ have been introduced in Accardi and Kozyrev \cite{AcKo01}. This paper motivated the general definitions of CPD-kernels and of CPD-semigroups in \cite{BBLS04}. (One should note that the definition in \cite{AcKo01} is slightly weaker, but this is compensated by that in \cite{AcKo01} the algebra is $\sB(H)$ and triviality of \nbd{W^*}correspondences over $\sB(H)$, where the weaker definition is equivalent to that in \cite{BBLS04}; see Skeide \cite[Lemma 5.27 and Remark 5.2.8]{Ske01} for an explanation.)

In \cite[Lemma 3.2.1]{BBLS04} there are listed some properties of a kernel $\eK$ that are equivalent to that $\eK$ is CPD. The most interesting are:
\begin{enumerate}
\item
For each finite choice of $\sigma_i\in S,a_i\in\cA,b_i\in\cB$ the map
\beqn{
a
~\longmapsto~
\sum_{i,j}b_i^*\eK^{\sigma_i,\sigma_j}(a_i^*aa_j)b_j
}\eeqn
is (completely) positive. (Apply positivity of these maps to $a=\U\ge0$ to see that $\eK$ is CPD. The nontrivial direction follows directly from Kolmogorov decomposition, because $a\mapsto\AB{x,ax}$ is clearly CP.)

\item
 For all choices $\sigma_1,\ldots,\sigma_n\in S$ $(n\in\N)$ the map
 \beqn{
 \eK^{(n)}
 \colon
\bfam{a_{ij}}
~\longmapsto~
\bfam{\eK^{\sigma_i,\sigma_j}(a_{ij})}
}\eeqn
from $M_n(\cA)$ to $M_n(\cB)$ is (completely) positive.
\end{enumerate}
For $\#S$ finite and fixed $n=\#S$, the second property with complete positivity is Heo's definition in \cite{Heo99}.
\enote
}

\section{Positivity, tensor product, and Schur product}\label{posSEC}

As we said in the introduction, a good possibility to check for positivity of something is to write that something as a square, that is, to find a square root. Looking at it from the other end, we can say a good notion of positivity is a notion that admits, as a theorem, the statement that every positive element has a whatsoever square root.

\nbd{C^*}Algebras do fulfill that criterion. An element in $b$ a \nbd{C^*}algebra $\cB$ is positive if and only if it can be written as $\beta^*\beta$ for some $\beta\in\cB$. An equivalent condition (frequently used as definition of positive element) is that $b=b^*$ and the spectrum $\sigma(b)$ is contained in $\R_+$ (\hl{spectral positivity}). Another equivalent condition is that  $\vp(b)\ge0$ for all positive linear functionals $\vp$ on $\cB$ (\hl{weak positivity}), where $\vp$ is \hl{positive} if $\vp(b^*b)\ge0$ for all $b\in\cB$.

Spectral positivity of $b\in\cB$ has the advantage that, apart from self-adjointness, it uses only the spectrum of $b$. And the spectrum of an element in a unital \nbd{C^*}algebra depends only the unital \nbd{C^*}subalgebra generated by that element. (The extra condition that a positive element must be self-adjoint, is something we gladly accept, looking at how powerful the other property is.) This means, no matter in how big another \nbd{C^*}algebra $\cA$ is into which we embed $\cB$, our element $b$, which is positive in $\cB$, continues being positive also in $\cA$.

Many other properties when dealing with positive elements are proved by spectral calculus. When constructing new Hilbert modules from given ones, it would be nice if checking on positivity of the new inner products in $\cB$, one would no longer need manipulations involving spectral calculus. Instead, one would simply refer to known results about positivity in $\cB$ or related \nbd{C^*}algebras and for the rest has purely algebraic operations. We illustrate what we mean, by giving several algebraic constructions that result in a completely algebraic proof of positivity of the inner product in the tensor product and positivity of the Schur product.

The \hl{direct sum} $E:=E_1\oplus\ldots\oplus E_n$ of Hilbert \nbd{\cB}modules with the inner product of $X=(x_1,\ldots,x_n)$ and $Y=(y_1,\ldots,y_n)$ defined as $\AB{X,Y}=\sum_{i}\AB{x_i,y_i}$ is a Hilbert \nbd{\cB}module. (Indeed, the sum of positive elements in a \nbd{C^*}algebra is positive. The rest follows as for Hilbert spaces. The direct sum of an infinite family of Hilbert modules still has to be completed.)

By $\sB^a(E,F)$ we denote the space of \hl{adjointable} operators from a Hilbert module $E$ to a Hilbert module $F$.  It is easy to check that
\beqn{
\sB^a(E)
~=~
\sMatrix{\sB^a(E_1,E_1)&\ldots&\sB^a(E_n,E_1)\\\vdots&&\vdots\\\sB^a(E_1,E_n)&\ldots&\sB^a(E_n,E_n)},
}\eeqn
acting in the obvious way on $X^n\in E$. Since $\sB^a(E)$ is a \nbd{C^*}algebra, for every $a\in\sB^a(E_i,E_j)$ the \it{square} $a^*a$ is a positive element of $\sB^a(E_i)$, because it is positive in $\sB^a(E)$.

\bex\label{linkex}
The \hl{linking algebra} $\tMatrix{\cB&E^*\\E&\sK(E)}\subset\sB^a(\cB\oplus E)$ for a Hilbert module over a \nbd{C^*}algebra $\cB$ is the most important example. Here $E^*=\CB{x^*\colon x\in E}\subset\sB^a(E,\cB)$ where $x^*\colon y\mapsto\AB{x,y}$, and the \nbd{C^*}subalgebra $\sK(E):=\cls EE^*$ of $\sB^a(E)$ is called the algebra of \hl{compact operators} on $E$.

It follows that $E^*$ is a correspondence from $\cB$ to $\sK(E)$ with inner product $\AB{x^*,y^*}=xy^*$. In particular, the element $xx^*$ is a positive element of $\sK(E)\subset\sB^a(E)$. Note that the algebra of compact operators on $E^*$ is $\sK(E^*)=\cls E^*E=\cls\AB{E,E}=:\cB_E$, the \hl{range ideal} of the inner product of $E$.
\eex

Now we have direct sums (in particular, we have the \hl{column space} $E^n$, the direct sum of $n$ copies of $E$), and we have the \hl{dual correspondence} $E^*$. So, nobody prevents us from combining these constructions. We define the \hl{row space} $E_n:=((E^*)^n)^*$. Observe that $E_n$ is a Hilbert module over the compact operators on $(E^*)^n$, that is, over $M_n(\sK(E^*)=M_n(\cB_E)$. The Hilbert module structure of the ideal $M_n(\cB_E)$ in $M_n(\cB)$, easily (and uniquely) extends to $M_n(\cB)$. After this, it is an easy exercise to verify that $E_n$ consists of elements $X_n=(x_1,\ldots,x_n)$ with inner product and module action given by
\baln{
\AB{X_n,Y_n}_{i,j}
&
~=~
\AB{x_i,y_i}
&
(X_nB)_i
&
~=~
\sum_jx_jb_{j,i}.
}\ealn
(Exercise: $\sB^a(E_n)=\sB^a(E)$ with $aX_n=(ax_1,\ldots,ax_n)$. In particular, $X_nY_n^*=\sum_ix_iy_i^*\in\sK(E)$.)

\bcor
For all choices of $n\in\N$ and $x_1,\ldots,x_n\in E$ the matrix $\bfam{\AB{x_i,x_j}}_{i,j}$ is a positive element of $M_n(\cB)$.
\ecor

Recall how this went. The inner product of direct sums (in particular, of column spaces) is positive, because the sum of positive elements in a \nbd{C*}algebra is positive. $xx^*$ is positive in $\sK(E)$ because it is positive in the linking algebra. Therefore, $E^*$ is Hilbert module over $\sK(E)$. The rest is purely algebraic manipulation, iterating this dualization operation with direct sums.

\brem
The corollary is a key ingredient for proving positivity of tensor products and Schur products. Of course, both $xx^*\ge0$ and the corollary can be proved in a different way, for instance, by making use of the well-known fact that $a\in\sB^a(E)$ is positive if and only if $\AB{x,ax}$ is positive for all $x\in E$; see \cite{Pas73}. But, the proof of this fact requires considerably more spectral calculus. Even in the scalar case $\cB=\C$, we think it will be difficult to find a simpler argument than ours above. But note that also in the scalar case we do have to recognize that the dual of a Hilbert space $H$, $H^*$, carries the structure of a Hilbert \nbd{\sK(H)}module.
\erem

To construct the \hl{tensor product} $E\odot F$ (over $\cB$) of a correspondence $E$ from $\cA$ to $\cB$ and a correspondence $F$ from $\cB$ to $\cC$ one tries to define a \nbd{\cC}valued semiinner product on $E\otimes F$ in the only reasonable way by
\beqn{
\AB{x\otimes y,x'\otimes y'}
~:=~
\AB{y,\AB{x,x'}y'}
}\eeqn
and sesquilinear extension. To show positivity we would be fine, if we had to show it only for elementary tensors. Indeed, the element $\AB{x,x}\in\cB$ is positive, so that we may write it as $\beta^*\beta$. Hence, $\AB{y,\AB{x,x}y}=\AB{y,\beta^*\beta y}=\AB{\beta y,\beta y}\ge0$. Now let $x_1,\ldots,x_n\in E$ and $y_1,\ldots,y_n\in F$ and put $X_n:=(x_1,\ldots,x_n)\in E_n$ and $Y^n:=(y_1,\ldots,y_n)\in F^n$. Observe that $F^n$ is a correspondence from $M_n(\cB)$ to $\cC$ with the obvious left action of $M_n(\cB)$. We find
\beqn{
\BAB{\sum_ix_i\otimes y_i,\sum_ix_i\otimes y_i}
~=~
\sum_{i,j}\AB{y_i,\AB{x_i,x_j}y_j}
~=~
\AB{X_n\otimes Y^n,X_n\otimes Y^n}
~\ge~
0.
}\eeqn
Once this is established, we may divide out the length-zero elements $\sN$, complete, and obtain $E\odot F$. We denote $x\odot y:=x\otimes y+\sN$. Note that $E\odot F$ is a correspondence from $\cA$ to $\cC$ which is determined up to bilinear unitary equivalence by the property that it is generated by elements $x\odot y$ having inner product $\AB{x\odot y,x'\odot y'}=\AB{y,\AB{x,x'}y'}$ and fulfilling $a(x\odot y)=(ax)\odot y$.

\bcor
The composition $\eL\circ\eK$ of CPD-kernels $\eK$ and $\eL$ as in Definition \ref{compdef}, is CPD, too.
\ecor

\proof
Denote by $(E,i)$ and $(F,j)$ the Kolmogorov decompositions of $\eK$ and $\eL$, respectively. Then
\beqn{
(\eL\circ\eK)^{\sigma,\sigma'}
~=~
\eL^{\sigma,\sigma'}\circ\eK^{\sigma,\sigma'}
~=~
\AB{j(\sigma),\AB{i(\sigma),\bullet i(\sigma')}j(\sigma')}
~=~
\AB{i(\sigma)\odot j(\sigma),\bullet i(\sigma')\odot j(\sigma')}
}\eeqn
is, clearly, a CPD-kernel.\qed

\bob
The preceding proof also shows that the GNS-correspondence of $\eL\circ\eK$ is the \nbd{\cA}\nbd{\cC}subcorrespondence of $E\odot F$ generated by all $i(\sigma)\odot j(\sigma)$, with $\sigma\mapsto i(\sigma)\odot j(\sigma)$ as cyclic map. Note that $E\odot F=(\cls\cA i(S)\cB)\odot(\cls\cB j(S)\cC)=\cls\,\bCB{\,a i(\sigma)\odot bj(\sigma')c\colon a\in\cA;b\in\cB;c\in\cC;\sigma,\sigma'\in S}$. So, $E\odot F$ is (usually much) bigger than the GNS-correspondence of $\eL\circ\eK$.
\eob

\bex\label{CPD-PSex}
A \hl{product system} is a family $E^\odot=\bfam{E_t}_{t\ge0}$ of \nbd{\cB}correspondences with an associative product $E_s\times E_t\ni(x_s,y_t)\mapsto x_sy_t\in E_{s+t}$ that extends as a bilinear unitary $E_s\odot E_t\rightarrow E_{s+t}$. If $\eT=\bfam{\eT_t}_{t\ge0}$ is a CPD-semigroup over $S$ on a unital \nbd{C^*}algebra $\cB$, then, by the preceding observation, the GNS-correspondences $\sE_t$ of the $\eT_t$ fulfill
\beqn{
(\sE_{s^n_{m_n}}\odot\ldots\odot\sE_{s^n_1})\odot\ldots\odot(\sE_{s^1_{m_1}}\odot\ldots\odot\sE_{s^1_1}).
~\supset~
\sE_{s^n_{m_n}+\ldots+s^n_1}\odot\ldots\odot\sE_{s^1_{m_1}+\ldots+s^1_1}
}\eeqn
For fixed $t>0$, this gives rise to an inductive limit over tuples $(t_n,\ldots,t_1)\in(0,\infty)^n$ with $t_n+\ldots+t_1=t$. For the resulting correspondences $E_t\supset\sE_t$ the inclusion $\sE_s\odot\sE_t\supset\sE_{s+t}$ becomes equality $E_s\odot E_t=E_{s+t}$. The elements $\xi^\sigma_t:=i_t(\sigma)\in\sE_t\subset E_t$ fulfill $\xi^\sigma_s\xi^\sigma_t=\xi^\sigma_{s+t}$, that is, for each $\sigma\in S$ the family ${\xi^\sigma}^\odot=\bfam{\xi^\sigma_t}_{t\ge0}$ is a \hl{unit}. Moreover, we have
\beqn{
\AB{\xi^\sigma_t,\bullet\xi^{\sigma'}_t}
~=~
\eT^{\sigma,\sigma'}_t
}\eeqn
for all $\sigma,\sigma'\in S$, and the set $\CB{{\xi^\sigma}^\odot\colon\sigma\in S}$ of units generates $E^\odot$ as a product system. We refer to $E^\odot$ as the \hl{GNS-system} of $\eT$ and to $\CB{{\xi^\sigma}^\odot\colon\sigma\in S}$ as the \hl{cyclic set of units}. We see:

\thmbox{
The square root of a CPD-semigroup (in particular, of a CP-semigroup) is a product system.
}
\eex

\SQ{\vspace{-3ex}
\bnote
For CP-semigroups (that is, for a one-point set $S$), the preceding construction is due to Bhat and Skeide \cite{BhSk00}. This seems to be the first publication where product systems of correspondences occur. The generalization to CPD-semigroup is from \cite{BBLS04}.

Meanwhile, structures like the family of GNS-correspondences $\sE_t$ with inclusions such as $\sE_{s+t}\subset\sE_s\odot\sE_t$ started to be investigated under the name of \hl{subproduct systems} by Shalit and Solel \cite{ShaSo09}, and under the name of \hl{inclusion systems} by Bhat and Mukherjee \cite{BhMu10}. In \cite{BhMu10}, which considers only subproduct systems of Hilbert spaces, it is proved (among other results) by the same inductive limit procedure that every subproduct system embeds into a proper product system. It is clear that this is true also for correspondences.
\enote
}
As a further application of the tensor product, we discuss a Stinespring construction for CPD-kernels.

\bex\label{Stinespring}
Let $\eK$ be CPD-kernel over a set $S$ from $\cA\ni\U_\cA$ to $\cB$, and denote by $(E,i)$ its Kolmogorov decomposition. Suppose $\cB$ is a concrete \nbd{C^*}algebra of operators acting nondegenerately on a Hilbert space $G$. In other words, suppose $G$ is a correspondence from $\cB$ to $\C$. Put $H:=E\odot G$. Observe that $H$ is a correspondence from $\cA$ to $\C$. In other words, $H$ is a Hilbert space with a nondegenerate action of $\cA$ that may be used to define the \hl{Stinespring representation} $\rho$ of $\cA$ on $H$ by $\rho(a)\colon x\odot g\mapsto a(x\odot g)=(ax)\odot g$. Note further that each $x\in E$ gives rise to an operator $L_x\colon g\mapsto x\odot g$. One readily verifies that
\beqn{
L_{i(\sigma)}^*\rho(a)L_{i(\sigma')}
~=~
\eK^{\sigma,\sigma'}(a)
}\eeqn
for all $\sigma,\sigma'\in S$ and $a\in\cA$.
\eex

\SQ{\vspace{-3ex}
\bnote\label{KSGNSnote}
Note that this easily generalizes to the case when we replace $\cB\subset\sB(G)$ with $\cB\subset\sB^a(F)$ where $F$ is some Hilbert \nbd{\cC}module. More precisely, $F$ is a correspondence from $\cB$ to $\cC$. In this case, also $H=E\odot F$ is a Hilbert \nbd{\cC}module, but for the rest nothing changes. For a one-point set $S$, this is known as \hl{KSGNS-construction}; see Lance \cite{Lan95} or Murphy \cite{Mur97}.
\enote
}
For a one-point set $S$ we recover the usual Stinespring construction for a CP-map from a unital \nbd{C^*}algebra into a concrete operator \nbd{C^*}algebra (cf.\ also Remark \ref{Stinerem}). However, neither the definition of CPD-kernel (or that of PD-kernel) nor the Kolmogorov decomposition for a CPD-kernel (or that of a PD-kernel) require that $\cB$ is represented as an operator algebra. CPD-kernels $\eK$ and $\eL$ (in particular, CP-maps) may be composed, and the Kolmogorov decomposition of the composed CPD-kernel $\eL\circ\eK$ can easily be recovered inside the tensor product $E\odot F$ of those $(E,i)$ and $(F,j)$ of the factors $\eK$ and $\eL$, respectively, as $_\cA i\odot j_\cC:=\cls\CB{ai(\sigma)\odot j(\sigma)c\colon a\in\cA,\sigma\in S,c\in\cC}$. On the contrary, \bf{never} will the Stinespring representation $\rho\colon a\mapsto a\odot\id_G$ of $\cA$ on $H=E\odot G$ for $\eK$ help to determine the Stinespring representation of $\cA$ for $\eL\circ\eK$!
\SQ{
Indeed, if $\cC\subset\sB(K)$, rather than the Stinespring representation of $\cA$ associated with the identity representation of $\cB$ on $G$, one would need the Stinespring representation of $\cA$ associated with the Stinespring representation $\pi\colon b\mapsto b\odot\id_K$ of $\cB$ on $L:=F\odot K$ for $\eL$. (By this we mean the representation $a\mapsto a\odot\id_L$ of $\cA$ on $E\odot L=E\odot F\odot K$.) The Stinespring representation of $\cA$ for $\eL\circ\eK$ would, then, be the representation $a\mapsto a\odot\id_L$ on $E\odot F\odot K$ restricted to the invariant subspace $({_\cA i\odot j_\cC})\odot K$. While this latter construction depends explicitly on $\eL$, or better on the Stinespring representation of $\cB$ associated with $\eL$, the GNS-correspondence of $\eK$ is universal and works for composition with all $\eL$.
}
Conclusion:
\thmbox{
Doing Stinespring representations for the individual members of a CP-semigroup on $\cB\subset\sB(G)$, is approximately as ingenious as considering a \nbd{2\times2}system of complex linear equations as a real \nbd{4\times4}system (ignoring all the structure hidden in the fact that certain \nbd{2\times2}submatrices are very special) and applying the Gauß algorithm to the \nbd{4\times4}system instead of trivially resolving the \nbd{2\times2}system by hand.
}

\bex
Let $E$ and $F$ be Hilbert modules over \nbd{C^*}algebras $\cB$ and $\cC$, respectively. Let $\vp$ be a map from $\cB$ to $\cC$. We say a linear map $T\colon E\rightarrow F$ is a \hl{\nbd{\vp}map} if
\beqn{
\AB{T(x),T(x')}
~=~
\vp(\AB{x,x'})
}\eeqn
for all $x,x'\in E$. Suppose $\cB$ is unital and $T$ is a \nbd{\vp}map for some CP-map $\vp$ from $\cB$ to $\cC$. Do the GNS-construction $(\sF,\zeta)$ for $\vp$. Then it is easy to verify that $v\colon x\odot b\zeta c\mapsto T(xb)c$ extends as an isometry $E\odot\sF\rightarrow F$. Of course,
\beqn{
v(x\odot\zeta)
~=~
T(x).
}\eeqn
Now let  $H_1,H_2$ denote Hilbert spaces. Put $\cC:=\sB(H_1)$, and put $F:=\sB(H_1,H_2)$. (This is a Hilbert \nbd{\sB(H_1)}module with inner product $\AB{y,y'}:=y^*y'$.) Assume $\Phi\colon E\rightarrow F=\sB(H_1,H_2)$ is a \nbd{\vp}map for the CP-map $\vp\colon\cB\rightarrow\cC=\sB(H_1)$, and construct the ingredients $\sF,\zeta,v$ as before.  Put $K_1:=\sF\odot H_1$, and denote by $\rho\colon b\mapsto b\odot\id_{H_1}$ the Stinespring representation of $\cB$ on $K_1$. Put $K_2:=E\odot K_1$, and denote by $\Psi\colon x\mapsto L_x=x\odot\id_{K_1}$ (the \hl{representation} of $E$ into $\sB(K_1,K_2)$ \hl{induced} by $\rho$). Denote $V:=L_\zeta=\zeta\odot\id_{H_1}\in\sB(H_1,K_1)$ and $W:=(v\odot\id_{H_1})^*\in\sB(F\odot H_1,E\odot\sF\odot H_1)=\sB(H_2,K_2)$. Then $K_1$; $K_2$; $V\in\sB(H_1,K_1)$; $W\in\sB(H_2,K_2)$; $\rho\colon\cB\rightarrow\sB(K_1)$; $\Psi\colon E\rightarrow\sB(K_1,K_2)$ fulfill the following:
\begin{enumerate}
\item
$W^*\Psi(x)V=T(x)$ for all $x\in E$.

\item
$\Psi(x)^*\Psi(x')=\rho(\AB{x,x'})$ for all $x,x'\in E$.

\item
$\rho$ is a nondegenerate representation.

\item
$W$ is a coisometry.
\end{enumerate}
A sextuple with these properties is determined uniquely up to suitable unitary equivalence.
\eex

\SQ{\vspace{-3ex}
\bnote
Existence and uniqueness of a sextuple fulfilling these properties is proved in Bhat, Ramesh, Sumesh \cite[Theorems 2.1 and 2.4]{BRS10p}. It should be noted that several additional conditions stated in \cite{BRS10p}, are automatic, once the four preceding conditions are fulfilled. (This result has been stated first by Asadi \cite{Asa09} with an unnecessary condition and with an incorrect proof.) The construction of $\sF,\zeta,v$ and the reduction of \cite{BRS10p} to it, is Skeide \cite{Ske10p1}.
\enote
}

\section{More Examples}

\bex
Every \nbd{\C}valued PD-kernel $\ek$ over a set $S$ gives rise to a CPD-kernel $\eK$ on $\sB(\C)=\C$ over $S$, if we interpret a complex number $z$ as a map $z\colon w\mapsto wz$ on $\C$. Under this identification the Kolmogorov decompositions of $\ek$ and of $\eK$ coincide. (After all, every Hilbert space is also a correspondence over $\C$ in the only possible way.) And the composition of the CPD-kernels corresponds to the usual Schur product of the PD-kernels. Of course, the tensor product of the Kolmogorov decompositions is the usual tensor product of Hilbert spaces.

In the same way, a \nbd{\cB}valued PD-kernel $\ek$ over $S$ can be interpreted as a CPD-kernel $\eK$ from $\C$ to $\cB$, by interpreting $b\in\cB$ as map $z\mapsto zb$.

About a CPD-kernel $\eK$ over $S$from $\cB$ to $\C$ we cannot say much more than that it is, by Kolmogorov decomposition, a representation of $\cB$ on a Hilbert space $G$ which is generated by $\cB$ from a set of $\#S$ vectors having certain inner products. Of course, for a one-point set $S=\CB{\om}$ we recover the GNS-construction for a positive linear functional $\vp\colon\cB\rightarrow\C$ with one vector $i(\om)$ that is \hl{cyclic} for $\cB$ (that is, $\ol{\cB i(\om)}=G$) such that $\AB{i(\om),\bullet i(\om)}=\vp$.
\eex

\bex
The tensor product of the Kolmogorov decompositions of (C)PD-kernels on $\C$ is the usual tensor product of Hilbert spaces. But the inclusion of the Kolmogorov decomposition for the composition into that tensor product remains. In particular, the construction of a product system, as discussed in Example \ref{CPD-PSex} for CP(D)-semigroups, works also here and the inductive limit is (almost always) proper. Just that now the result is a product system of Hilbert spaces (\hl{Arveson system}, henceforth).

It is noteworthy that under very mild conditions (just measurability of the semigroups $\et^{\sigma,\sigma'}$ in $\C$) the Arveson system of a (C)PD-semigroup (or Schur semigroup of positive definite kernels) $\et$ can be computed explicitly. Indeed, under this condition we may define the generator $\el$ of $\et$ as the kernel
\beqn{
\el^{\sigma,\sigma'}
~:=~
\textstyle\frac{d}{dt}\big|_{t=0}\et^{\sigma,\sigma'}_t.
}\eeqn
One easily verifies that $\el$ is \hl{conditionally positive definite}, that is,
\beqn{
\sum_{i,j}\bar{z}_i\el^{\sigma_i,\sigma_j}z_j
~\ge~
0
}\eeqn
whenever $\sum_iz_i=0$. Note, too, that $\el$ is \hl{hermitian} in the sense that $\ol{\el^{\sigma,\sigma'}}=\el^{\sigma',\sigma}$ for all $\sigma,\sigma'\in S$ (simply because $\et_t$, obviously, is hermitian).

Note that for every choice of $\beta_\sigma\in\C$ $(\sigma\in S)$ with $\et$, also $\tilde{\et}$ defined by setting $\tilde{\et}^{\sigma,\sigma'}_t:=e^{\ol{\beta}_\sigma t}~\et^{\sigma,\sigma'}_t~e^{\beta_{\sigma'}t}$ is a PD-semigroup with generator $\tilde{\el}^{\sigma,\sigma'}=\el^{\sigma,\sigma'}+\beta_{\sigma'}+\ol{\beta}_\sigma$. Fix a $\sigma_0\in S$. Then choose $\beta_{\sigma_0}$ such that $\tilde{\et}^{\sigma_0,\sigma_0}_t=1$, that is, such that $\tilde{\el}^{\sigma_0,\sigma_0}=\el^{\sigma_0,\sigma_0}+\beta_{\sigma_0}+\ol{\beta}_{\sigma_0}=0$. Further, for that $\beta_{\sigma_0}$ choose the other $\beta_\sigma$ such that $\tilde{\et}^{\sigma_0,\sigma}_t=1$, that is, such that $\tilde{\el}^{\sigma_0,\sigma}=\el^{\sigma_0,\sigma}+\beta_{\sigma_0}+\ol{\beta}_\sigma=0$. Verify that the kernel $\tilde{\el}$ we obtain in that way is not only conditionally positive definite, but really positive definite. Do the Kolmogorov decomposition $(K,i)$ for that kernel. Define $E_t:=\G(L^2(\SB{0,t},K))$ (symmetric Fock space) and observe that the $E_t$ from an Arveson system via 
\beqn{
E_s\otimes E_t
~\rightarrow~
\s_tE_s\otimes E_t
~=~
E_{s+t},
}\eeqn
where the first step is the second quantized time-shift emerging from $\SB{0,s}\mapsto\SB{t,t+s}$, and where the second step is the usual factorization $\G(H_1)\otimes\G(H_2)=\G(H_1\oplus H_2)$. Note that for $\beta\in\C$ and $k\in K$ the elements $\xi_t(\beta,k):=e^{\beta t}\psi(\I_{\SB{0,t}}k)\in E_t$ ($\psi$ denoting exponential vectors and $\I$ indicator functions) form units. Finally,
\beqn{
\AB{\xi_t(-\beta_\sigma,i(\sigma)),\xi_t(-\beta_{\sigma'},i(\sigma'))}
~=~
e^{-\ol{\beta}_\sigma t}e^{\tilde{\el}^{\sigma,\sigma'}t}e^{-\beta_{\sigma'}t}
~=~
\et^{\sigma,\sigma'}_t.
}\eeqn
Since the $i(\sigma)$ generate $K$ and since $i(\sigma_0)=0$, these units generate the whole product system; see Skeide \cite{Ske00a}. In other words, we have just constructed the Kolmogorov decomposition of $\et$. Moreover, a brief look at the construction shows, that it actually does not depend on that $\el$ is \it{a priori} the generator of a PD-semigroup, but only on its properties to be hermitian and conditionally positive definite. We, thus, also have proved that every such kernel generates a PD-semigroup by exponentiation. This relation is called the \hl{Schönberg correspondence} between PD-semigroups and conditionally positive definite hermitian \nbd{\C}valued kernels.
\eex

\SQ{\vspace{-3ex}
\bnote
Without notions like product systems and units for them (these came not before Arveson \cite{Arv89}), the possibility to realize PD-semigroups as inner products of suitably normalized exponential vectors has been discovered as early as Parthasarathy and Schmidt \cite{PaSchm72}. They applied it to characteristic functions (the Fourier transform) of the convolution semigroup of distributions of Lévy processes, and used it to prove representability of an arbitrary Lévy process (starting at $0$) by (possibly infinite, for instance, for the Cauchy process) linear combinations of the usual creation, conservation, and annihilation processes on the Fock space; a result that Schürmann \cite{MSchue93} generalized to quantum Lévy processes.
\enote

\bnote
A natural question is if the a similar representation result holds for CPD-semi\-groups $\eT$ on $\cB$, when we replace the symmetric Fock space by the time ordered Fock module; \cite{BhSk00}. The answer is: Yes, for von Neumann algebras, and if the CPD-semigroup is uniformly continuous; \cite{BBLS04}. The procedure is essentially the same. One has to find a $\sigma_0$ such that the semigroup $\eT^{\sigma_0,\sigma_0}$ can be normalized to give the trivial semigroup $\id$ on $\cB$. However, as this is no longer possible in general, one rather has to try add to a point $\sigma_0$ to $S$ and to extend the CPD-semigroup over $S$ to a CPD-semigroup over $S\cup\CB{\sigma_0}$ in a consistent way. In the von Neumann case this is always possible. But unlike the scalar case, here this result is a hard to prove. In fact, it is equivalent to the deep result by Christensen and Evans \cite{ChrEv79} who found the form of the generator of a uniformly continuous CP-semigroup on a von Neumann algebra, and in \cite{BBLS04} only equivalence to \cite{ChrEv79} is proved. From this point on, everything goes like the scalar case. The unit $\om^\odot$ representing the value $\sigma_0$, actually already belongs to the GNS-system of $\eT$, and what is generated by all these units is a whole product system of time ordered Fock modules.

In the \nbd{C^*}case the situation is worse. Not only is it possible that the GNS-system embeds into a product system of time ordered Fock modules, but is not isomorphic to one. It is even possible that a the GNS-system of a uniformly continuous CP-semigroup $T$ (even an automorphism semigroup) does not embed at all into a product system of time ordered Fock modules. In fact, a the GNS-system embeds if and only if the CP-semigroup is \hl{spatial}. Spatial means that $T$ dominates a CP-semigroup of the form $b\mapsto c_t^*bc_t$ for a norm continuous semigroup of elements $c_t$ in $\cB$. All this has been discussed for CP-semigroups in Bhat, Liebscher, and Skeide \cite{BLS10}. The extension of the discussion to CPD-semigroup can be found in Skeide \cite{Ske10}.

It should be noted that in Skeide \cite{Ske06d} (preprint 2001) we called a product system $E^\odot$ \hl{spatial} if it admits a unit $\om^\odot$ that gives the trivial semigroup $\AB{\om_t,\bullet\om_t}=\id_\cB$ on $\cB$ (like the unit ${\xi^{\sigma_0}}^\odot$ corresponding to the index $\sigma_0$ of the extended CPD-semigroup above). Spatial product system allow a classification scheme most similar to Arveson systems. They are type I (that is, Fock) and type II (non-Fock) depending on whether or not they are generated by a set of units that determines a uniformly continuous CPD-semigroup. And they have an index (the one-particle sector of the maximal Fock part) that behaves additive (direct sum) under a suitable spatial product of spatial product systems; see \cite{Ske06d}. For Arveson systems (which, unlike general product systems, possess a tensor product) it may but need not coincide with the tensor product of Arveson systems,
\enote
}

\bex\label{1dimex}
Of course, also a semigroup $\vt=\bfam{\vt_t}_{t\ge0}$ of unital endomorphisms of $\cB$ is a CP-semigroup. One readily verifies that its GNS-system consists of the correspondences $\cB_t=\cB$ as right Hilbert \nbd{\cB}module, but with left action $b.x_t:=\vt_t(b)x_t$. The generating unit is $\U_t=\U$. (Note that each $\cB_t$ is the GNS-correspondence of $\vt_t$. No inductive limit is necessary, because $\cB_s\odot\cB_t=\cB_{s+t}$.) On the other hand, if we have such a \hl{one-dimensional} product system $\cB_t$, then $\vt_t(b):=b.\U_t$ defines a unital endomorphism semigroup.
\eex

\SQ{\vspace{-3ex}
\bnote
This might appear to be a not so interesting product system. However, it turns out that every product system of correspondences over $\cB$ arises from such a one-dimensional product system of correspondences over $\sB^a(E)$. (Anyway, one should be alarmed by the fact that a one-dimensional product system of correspondences over $\cB$ is isomorphic to the trivial one if and only if the endomorphism semigroup consists of inner automorphisms of $\cB$.) The relation between these two product systems, is an operation of Morita equivelance; see the following Exampe \ref{BaEex} where the operation is scratched, and Note \ref{BaEnote} where Morita equivalence is mentioned very briefly. This scratches the intimate relationship between \nbd{E_0}semigroups (that is, unital endomorphism semigroups) on $\sB^a(E)$ and product systems (first Arveson \cite{Arv89,Arv90} for \nbd{E_0}semi\-groups on $\sB(H)$ and Arveson systems), and we have no space to give an account. Instead of a long list of papers, we mention Skeide \cite{Ske08p1} where the theory has been completed, and the references therein.
\enote
}

\bex\label{BaEex}
In Example \ref{Stinespring} and the remarks following it, we have emphatically underlined that we do not consider it too much of a good idea to tensor a GNS-correspondence (from $\cA$ to $\cB$, say) with a representation space $G$ (or module $F$) of $\cB$. This changes if actually $\cB=\sB(G)$ or $\cB=\sB^a(F)$.

To fix letters a bit more consistently, let $E$, $F$, and $G$ be a Hilbert \nbd{\cB}module, a Hilbert \nbd{\cC}module, and a Hilbert \nbd{\cD}module, respectively. And suppose $\sE$ and $\sF$ are correspondences from $\sB^a(E)$ to $\sB^a(F)$, and from $\sB^a(F)$ to $\sB^a(G)$, respectively. If $\sE$ is the GNS-correspondence of a CP-map (or even of a unital homomorphism), then nobody would be surprised that we require technical conditions for that CP-map (like normality in the case when $E$ is a Hilbert space) that are reflected also by the left action of that correspondence. In the the framework of \nbd{C^*}modules this condition is \hl{strictness} of the left action of $\sB^a(E)$ (or, more precisely, strictness on bounded subsets). It is equivalent to that the compacts $\sK(E)$ alone act already nondegenerately on $\sE$. So, let us also suppose that $\sE$ and $\sF$ are strict correspondences in that sense. Let us compute $E^*\odot\sE\odot F$ and $F^*\odot\sF\odot G$. (Recall the definition of the dual correspondence $E^*$ from Example \ref{linkex}.) Then
\bmun{
(E^*\odot\sE\odot F)\odot(F^*\odot\sF\odot G)
~=~
E^*\odot\sE\odot(F\odot F^*)\odot\sF\odot G
\\
~=~
E^*\odot\sE\odot\sK(F)\odot\sF\odot G
~=~
E^*\odot(\sE\odot\sF)\odot G.
}\emun
Here $F\odot F^*=\sK(F)$ via $y'\odot y^*\mapsto y'y^*$, and in the last step we made use of strictness of the left action of $\sF$.

We see that ``tensor-sandwiching'' from the left and right with the respective representation modules is an operation that respects tensor products (and is, obviously, as associative as one may wish). Suppose $\sE^\odot$ is a product system of strict correspondences over $\sB^a(E)$. Then the correspondences $E_t:=E^*\odot\sE_t\odot E$ form a product system $E^\odot$ of correspondences over $\cB$. If $\sE^\odot$ is the one-dimensional product system of a strict \nbd{E_0}semigroup on $\sB^a(E)$ as in Example \ref{1dimex}, then $E^\odot$ is indeed the product system of \nbd{\cB}correspondences associated with that \nbd{E_0}semigroup. In the same way one obtains the product system of \nbd{\cB}correspondences of a strict CP(D)-semigroup group on $\sB^a(E)$ from the GNS-system of that CP(D)-semigroup.

It is noteworthy that by this elegant and simple method, we recover for $\sB(H)$ the constructions of Arveson systems from \nbd{E_0}semigroups on $\sB(H)$ or from CP-semigroups on $\sB(H)$. However, even in that simple case it is indispensable to understand the dual correspondence of the Hilbert space $H$, $H^*$, as a fully qualified Hilbert \nbd{\sB(H)}module.
\eex

\SQ{\vspace{-3ex}
\bnote\label{BaEnote}
For \nbd{E_0}semigroups of $\sB(H)$ this product system has been constructed in a different way by Bhat \cite{Bha96}. It is anti-isomorphic to the Arveson system from \cite{Arv89} (and need not be isomorphic to it; see Tsirelson \cite{Tsi00p1}). We imitated Bhat's construction in Skeide \cite{Ske02}; but since this construction requires existence of a unit vector in $E$, it is not completely general. The general construction above is from Skeide \cite{Ske09} (preprint 2004). The construction of the Arveson system of a CP-semigroup on $\sB(H)$ is again due to \cite{Bha96}; it goes via dilation of the CP-semigroup to an \nbd{E_0}semigroup and, then, determining the Arveson system of that \nbd{E_0}semigroup. In Skeide \cite{Ske03c} we gave a direct construction of that product system along the above lines. Bhat, Liebscher, and Skeide \cite{BLS08} discuss the generalization from \cite{Ske03c} to CP-semigroups on $\sB^a(E)$, Skeide \cite{Ske10} the case of CPD-semigroups.

``Tensor-sandwiching'', as we called it above, is actually \it{cum grano salis} an operation of Morita equivalence for correspondences. (For von Neumann algebras and modules, it is Morita equivalence. For \nbd{C^*}algebras and modules it is Morita equivalence up to strict completion.) This has been explained in \cite{Ske09} and fully exploited in \cite{Ske08p1} to complete the theory of classification of \nbd{E_0}semigroups by product systems.
\enote
}

\section{Positivity in $*$--algebras}

We said that if a notion of positivity is good, then positive things should have a square root. Positive elements of a \nbd{C^*}algebra have a square root inside the \nbd{C^*}algebra. In Skeide \cite{Ske01} we worked successfully with the definition that an element in a pre-\nbd{C^*}algebra $\cB$ is positive if it is positive in the completion of $\cB$, $\ol{\cB}$. In other words, such an element has a square root not necessarily in $\cB$, but it always has one in $\ol{\cB}$.

\bex
If polynomials behave nicely with respect to positivity, depends on where we wish to evaluate them. If $p$ is a polynomial with complex coefficients and we consider functions $x\mapsto p(x)$ from $\R$ to $\C$, then $p$ is positive in the sense that $p(x)\ge0$ for all $x\in\R$ if and only if there is a polynomial $q$ such that $p=\ol{q}q$. If we evaluate in $\C$ instead of $\R$ (requiring $p(z)\ge0$ for all $z\in\C$), then Liouville's theorem tells us that $p$ is constant. If we evaluate in a bounded interval $I=\SB{a,b}$ only, then we get more positive polynomials. For instance, the polynomial $p(x)=-x$ is positive on the interval $\SB{-1,0}$. But no $q$ will ever give back $p$ as a square $\ol{q}q$. Note that such polynomials that are positive on bounded intervals occur (and are necessary!) in typical proofs of the spectral theorem, where they serve to approximate indicator functions of intervals from above. Of course, $-x$ has a square root in the \nbd{C^*}algebra $C\SB{-1,0}$.
\eex

We already showed positivity of an element in $\sB^a(E)$ by writing it as $a^*a$ for some $a\in\sB^a(E,F)$. So, it sometimes is convenient even to leave completely the algebra under consideration. Our thesis was also that  Hilbert modules can be good square roots of positive-type maps.

It is our scope to propose here a new and flexible definition of positivity, that complements the algebraic notion of positivity from Accardi and Skeide \cite{AcSk08} (applied successfully in \cite{AcSk00a} to the ``square of white noise''). We will show that this notion allows for nice square roots of positive things.

For the balance of this section, $\sL(V,W)$ stands for \hl{linear} maps between vector spaces, and $\sL^a(E,F)$ stands for \hl{adjointable} maps between pre Hilbert \nbd{\cB}modules. The latter are linear automatically, but in general not necessarily bounded (unless at least one, $E$ or $F$, is complete).

\bdefi\label{SWposdef}
Let $\cB$ be a unital \nbd{*}algebra and let $\cS$ be a set of positive linear functionals on $\cB$. We say $b\in\cB$ is \hl{\nbd{\cS}positive} if $\vp(c^*bc)\ge0$ for all $\vp\in\cS,c\in\cB$. We say $\cB$ is \hl{\nbd{\cS}separated} if the functionals in $\vp$ separate the points of $\cB$ in the sense that $\vp(c^*bc')=0$ for all $\vp\in\cS$ and all $c,c'\in\cB$ implies $b=0$.
\edefi

\bob\label{posob}
By polarization and since every $b\in\cB$ is a linear combination of two self-adjoint elements, we may check separation only for $\vp(c^*bc)$ and $b=b^*$.
\eob

\bob
$\vp(c^*bc')=0$ for all $c,c'\in\cB$ means that the GNS-representation of the positive functional $\vp$ sends $b$ to $0$. In other words, if $\cB$ is \nbd{\cS}separated, then the direct sum of all GNS-representations for $\vp\in\cS$ is a faithful representation. That is, we may and will interpret $\cB$ as a concrete \nbd{*}algebra of operators $\cB\subset\sL^a(G)$ where the pre-Hilbert space $G$ is the (algebraic, of course) direct sum over the GNS-representations for all $\vp$. By Observation \ref{posob}, an element on $b\in\cB$ is \nbd{\cS}positive if and only if $\AB{g,bg}\ge0$ for all $g\in G$.

One might be tempted to try the latter as a definition of positivity for an arbitrary \nbd{*}algebra of operators. However, our representation has something special, namely, it is the direct sum of cyclic representations. It it this property that will allow us to work.
\eob

\brem
If $\cB$ is not \nbd{\cS}separated, then we simply may quotient out the kernel of the representation on $G$. After all, $\cS$ will typically contain those states that correspond to all possible measurements that the system described by our algebra $\cB$ allows. If these measurements do not separate two points $b_1$ and $b_2$ from $\cB$, it is pointless to consider them different.
\erem

In the sequel, we fix a unital \nbd{*}algebra $\cB$ that is \nbd{\cS}separated by a set $\cS$ of positive functionals and represented faithfully as operator algebra $\cB\subset\sL^a(G)$ on the direct sum $G$ of the GNS-spaces of all $\vp\in\cS$. 

The following Kolmogorov decomposition is the generalization of Example \ref{Stinespring} to the algebraic situation of this section.

\bthm\label{SWposthm}
Let $\cA$ be a unital \nbd{*}algebra. For some set $S$ let $\eK\colon S\times S\rightarrow\sL(\cA,\cB)$ be a CPD-kernel over $S$ from $\cA$ to $\cB$, in the sense that all sums in \eqref{CPD} are \nbd{\cS}positive. Then there exists a pre-Hilbert space $H$ with a left action of $\cA$, and a map $i\colon S\rightarrow\sL^a(G,H)$ such that
\beqn{
\eK^{\sigma,\sigma'}(a)
~=~
i(\sigma)^*ai(\sigma')
}\eeqn
for all $\sigma,\sigma'\in S$ and $a\in\cA$. We refer to $E:=\ls\cA i(S)\cB\subset\sL^a(G,H)$ as the \hl{GNS-correspondence} of $\eK$ with \hl{cyclic map} $i$, and to $(E,i)$ its \hl{Kolmogorov decomposition}.
\ethm

\proof
On $\cA\otimes S_\C\otimes\cB$ we define a sesquilinear map precisely as in the proof of Theorem \ref{KdABthm}. The map is right \nbd{\cB}linear, and the left action of $\cA$ is a \nbd{*}map. It is rather a tautology that $\AB{\bullet,\bullet}$ is positive for whatever notion of positivity we define on $\cB$. So, we may call that map a \hl{semiinner product} and $\cA\otimes S_\C\otimes\cB$ a \hl{semicorrespondence} from $\cA$ to $\cB$.

The true work starts now when we wish to divide out the length-zero elements $\sN:=\CB{x\in\cA\otimes S_\C\otimes\cB\colon\AB{x,x}=0}$. Obviously, $x\in\sN$ implies $xb\in\sN$. Next we show that
\beq{\label{CSIc}
x\in\sN
~\Longrightarrow~
\AB{x,y}=0
\text{~~~for all~~~}y\in\cA\otimes S_\C\otimes\cB.
}\eeq
Indeed, $\abs{\vp(c^*\AB{x,y}c')}^2=\abs{\vp(\AB{xc,yc'})}^2\le\vp(\AB{xc,xc})\vp(\AB{yc',yc'})=0$, so $\AB{x,y}=0$. From this, we immediately conclude that $x,y\in\sN$ implies $x+y\in\sN$, so $\sN$ is a right \nbd{\cB}submodule. Suppose again $x\in\sN$. Then $\AB{ax,ax}=\AB{x,a^*ax}=0$, by \eqref{CSIc}, so $\sN$ is an \nbd{\cA}\nbd{\cB}submodule. We, therefore, may define the \nbd{\cA}\nbd{\cB}module $E:=(\cA\otimes S_\C\otimes\cB)/\sN$. Now, if $x,y\in\cA\otimes S_\C\otimes\cB$ and $n,m\in\sN$. Then $\AB{x+n,y+m}=\AB{x,y}$, again by \eqref{CSIc}. Therefore, by $\AB{x+\sN,y+\sN}:=\AB{x,y}$ we well-define a semiinner product on $E$ which is \hl{inner} (that is $\AB{x,x}=0$ implies $x=0$ for all $x\in E$). We, thus, may call $E$ a \hl{precorrespondence} from $\cA$ to $\cB$. Of course, the elements $j(\sigma):=\U\otimes e_\sigma\otimes \U+\sN$ fulfill $\AB{j(\sigma),aj(\sigma')}=\eK^{\sigma,\sigma'}(a)$.

Next we wish to construct $E\odot G$. To that goal, we define as usual a sesquilinear map on $E\otimes G$ by setting $\AB{x\otimes g,x'\otimes g'}:=\AB{g,\AB{x,x'}g'}$. We have to face the problem to show that
\beqn{
\sum_{i,j}\AB{x_i\otimes g_i,x_j\otimes g_j}
~=~
\sum_{i,j}\AB{g_i,\AB{x_i,x_j}g_j}
~\ge~
0
}\eeqn
for all finite choices of $g_i\in G$ and $x_i\in E$. Recall that $G$ decomposes into a direct sum of pre-Hilbert subspaces $G_\vp$ with cyclic vectors $g_\vp$, say, such that $\AB{g_\vp,bg_{\vp'}}=0$ for all $\vp\ne\vp'$ and such that $g=\sum_\vp b_\vp g_\vp$ for suitable $b_\vp$ (different from $0$ only for finitely many $\vp$). It follows that
\beqn{
\sum_{i,j}\AB{g_i,\AB{x_i,x_j}g_j}
~=~
\sum_\vp\sum_{i,j}\AB{b^i_\vp g_\vp,\AB{x_i,x_j}b^j_\vp g_\vp}
~=~
\sum_\vp\BAB{ g_\vp,\bAB{\textstyle\sum_ix_ib^i_\vp,\sum_ix_ib^i_\vp}g_\vp}
~\ge~
0.
}\eeqn
Now it is clear that every $x\in E$ defines an operator $L_x\colon g\mapsto x\odot g$ with adjoint $L_x^*\colon y\odot g\mapsto\AB{x,y}g$. We end the proof, by putting $i(\sigma):=L_{j(\sigma)}$.\qed
\SQ{\vspace{-3ex}
\bnote\label{CSIrem}
Observe that \eqref{CSIc} functions as substitute for Cauchy-Schwarz inequality. We wish to underline that at least as important as the estimates that follow from Cauchy-Schwarz inequality for pre-Hilbert modules (for instance, that the norm of a Hilbert modules is a norm, or that the operator norm of $\sB^a(E)$ is a \nbd{C^*}norm), is that Cauchy-Schwarz inequality is true also for semiinner products. Note, too, that unlikefor Hilbert spaces, for semi-Hilbert modules the case $\AB{x,x}=0=\AB{y,y}$ in the proof Cauchy-Schwarz inequality require a small amount of additional work; see, for instance, \cite[Proposition 1.2.1]{Ske01}.
\enote
}

\vspace{1ex}\noindent
We should note that when wish to compose such kernels, so that also the  algebra to the left needs a positivity structure, then the CPD-condition must be supplemented with a compatibility condition for that positivity. (Compare this with the solution for the algebraic definition of positivity in \cite{AcSk08}.) We leave this to future work and close with a corollary about existence of square roots for \nbd{\cS}positive elements of $\cB$.

\bcor
\begin{enumerate}
\item
Suppose that $\cA=\C$ and $S=\om$ and choose an \nbd{\cS}positive $b\in\cB$. Then $\beta=i(\om)$ is an adjointable operator $G\rightarrow H$ such that $b=\beta^*\beta$.

\item
By Friedrichs' theorem and spectral calculus for self-adjoint operators, there exists a positive self-adjoint operator $\sqrt{b}\colon\ol{G}\supset D\rightarrow\ol{G}$ with $G\subset D$. Then $\beta:=\sqrt{b}\upharpoonright G$ and $H:=\sqrt{b}G$ is a possible choice for (1).
\end{enumerate}
\ecor

\proof
Only for (2) there is something to show. Indeed, after having chosen $\sqrt{b}$ and defined $\beta$ and $H$ as stated, let us choose some $\beta'$ and $H'$ that exist according to (1). If necessary replace $H'$ with $\beta'G$, so that $\beta'$ is surjective. Then $v\colon\beta'g\mapsto\sqrt{b}g=\beta g$ defines an isometry into $\ol{G}$ and a unitary onto $H$. Since $\beta'$ has an adjoint $H'\rightarrow G$, it follows that $\beta\colon G\rightarrow H$ has the adjoint $\beta^*={\beta'}^*v^*$. Of course, $\beta^*\beta g={\beta'}^*v^*\beta g={\beta'}^*\beta'g=bg$.\qed

\setlength{\baselineskip}{2.5ex}

\newcommand{\Swap}[2]{#2#1}\newcommand{\Sort}[1]{}
\providecommand{\bysame}{\leavevmode\hbox to3em{\hrulefill}\thinspace}
\providecommand{\MR}{\relax\ifhmode\unskip\space\fi MR }
\providecommand{\MRhref}[2]{%
  \href{http://www.ams.org/mathscinet-getitem?mr=#1}{#2}
}
\providecommand{\href}[2]{#2}

\lf\noindent
Michael Skeide: \it{Dipartimento S.E.G.e S., Università degli Studi del Molise, Via de Sanctis, 86100 Campobasso, Italy},
E-mail: \href{mailto:skeide@unimol.it}{\tt{skeide@unimol.it}},\\
Homepage: \url{http://www.math.tu-cottbus.de/INSTITUT/lswas/_skeide.html}

\end{document}